\documentclass[twoside,draft,reqno]{birkart}
\usepackage{amssymb}
\def\BC{{\bf C}}
\def\BP{{\bf P}}
\def\BR{{\bf R}}
\def\BZ{{\bf Z}}
\def\CM{{ M}}
\def\CN{{ N}}
\def\CO{{ O}}
\def\CW{{ W}}
\def\II{{II}}
\def\WP{{WP}}
\def\grad{{\rm grad\ }}
\def\Tr{{\rm Tr\ }}
\def\Im{{\rm Im\ }}
\def\Ext{{\rm Ext}}
\def\Hom{{\rm Hom}}
\begin{document}
\setcounter{page}{1}
\title[D-Branes on CY Manifolds]
{D-Branes on Calabi-Yau Manifolds}
\author[M.~R.~Douglas]{Michael~R.~Douglas}
\address{
Department of Physics and Astronomy,\br
Rutgers University,\br
Piscataway, NJ 08855--0849\br
and\br
I.H.E.S.,\br
Le Bois-Marie, Bures-sur-Yvette, 91440 France}
\email{mrd@physics.rutgers.edu}
\begin{abstract}
We give an overview of recent work on Dirichlet branes on Calabi-Yau
threefolds which makes contact with Kontsevich's homological mirror
symmetry proposal, proposes a new definition of stability which is
appropriate in string theory, and provides concrete quiver categories
equivalent to certain categories of branes on CY.
\end{abstract}
\maketitle
\section{Introduction}
Dirichlet branes, discovered by Dai, Leigh and Polchinski in 1989,
have been a central part of the dramatic progress in superstring theory
of recent years.  They are the simplest solitons in string theory;
the basic definition of Dirichlet brane (D-brane) 
is that it is an allowed boundary condition for a string. \cite{Pol}

Although they play many physical roles, from a mathematical point of
view their most salient feature is that they provide the natural
context in which Yang-Mills theory is embedded in string theory.
This has allowed physicists to rederive and better understand some
of the most beautiful mathematical constructions which have been
discovered in this area, such as the ADHM construction of instantons and
Nahm's construction of monopoles.  

D-branes on Calabi-Yau (CY) 
manifolds have been the focus of a number of recent
works (see \cite{Doug2,DD} for additional references).
Early on it was realized that this theory had a close
connection with mirror symmetry and especially 
Kontsevich's homological mirror
symmetry proposal \cite{Kont}, 
which interprets mirror symmetry as an equivalence between
the derived category of coherent sheaves (naturally associated to
Yang-Mills on K\"ahler manifolds) and a derived category proposed
by Fukaya, associated to isotopy classes of Lagrangian submanifolds.
These two classes of objects (usually called ``B'' and ``A''-type
branes respectively) correspond
to the two classes of supersymmetry preserving (or BPS) branes on a CY,
and thus the physical understanding of mirror symmetry indeed should imply
such an equivalence.

In physical terms, Kontsevich's proposal addresses only the structure
of topologically twisted open string theory (we explain this below),
while one expects more structure to appear in discussing the full
open string theory which governs the physics of D-branes.

One approach to this extra structure is to study the additional conditions
on the D-brane world-volumes which follow from space-time supersymmetry.
In the large volume limit, these two conditions are the hermitian Yang-Mills
equations (for B branes) and the special Lagrangian condition (for A branes).
In this context, Strominger, Yau and Zaslow formulated mirror symmetry as
the statement that the mirror $\CW$ to a CY $\CM$
could be obtained as the moduli space of a
particular A brane mirror to the B brane which is the point on $\CM$. 
\cite{SYZ}
This brane is expected to come in families which give a $T^3$ fibration
structure to $\CM$, a claim which has been proven in many cases by
Gross and by Ruan. \cite{Zhar,Gross,Ruan,Morr}

Although this approach clearly captures some essential part of the truth,
it should be realized that the specific additional conditions which were
used, the hermitian Yang-Mills and special Lagrangian conditions, are not
believed to be the correct physical conditions except in the large volume
limit.  This is related to the point that in string theory, CY manifolds
are not really Ricci-flat manifolds except in the large volume limit; the
Einstein equations are known to be corrected.  Indeed, the 
stringy version of CY K\"ahler moduli space, which is best thought of
as the complex structure moduli space of the mirror,
looks quite different from the conventional geometric version.

The correct replacements for either Ricci-flatness or the
supersymmetry conditions on D-branes are not known at this writing,
but some things are known, which we summarize here.  In particular,
one has strong reasons to believe that the metric defined by the SYZ
construction (which is the ``D$0$-metric'' of \cite{Doug1}
on the mirror manifold) is {\it not} Ricci-flat.

Having cast doubt on the usual mathematical starting point for the
discussion of D-branes, we are obliged to say something about what
replaces it.  Now the fundamental physical definition is
quite clear -- a D-brane is a boundary condition in the world-sheet
conformal field theory (CFT); the A and B branes each correspond to boundary
conditions preserving a particular half of world-sheet supersymmetry.

This definition has some conceptual advantages; notably, as in the
discussion of Greene and Plesser \cite{GP}, mirror symmetry is manifest.
The problem with it is that it is rather hard to make precise contact
with geometry except in the large volume limit.  However, a number of
general results have been obtained, which we outline.  In keeping with
the CFT framework, we freely use mirror symmetry in our considerations,
while regarding the precise connection with geometry as following from
subsequent considerations which are not yet completely understood.

An example of a geometrical statement which is on a sound footing in
string theory is the statement that the holomorphic properties of
B branes are independent of K\"ahler moduli, and thus are ``the same''
as in the large volume limit \cite{BDLR,Doug1}.  The precise meaning of
``same'' turns out to be subtle however and this is the point where
the derived category must enter the discussion.
A recent development is a precise contact with the graded categories
which appear in Kontsevich's proposal, which shows that (if we consider B
branes), as the K\"ahler moduli are varied, the gradings undergo a 
``flow'' determined by the variation of certain phases associated with
the branes (in A brane language, the phase appearing in the special
Lagrangian condition). \cite{DFR2,Doug3}

This puts us in a position to explain perhaps the simplest physical
structure carried by D-branes which is not already implicit in
\cite{Kont}.  This is the concept of ``marginal stability'' which
governs the possible variation of the spectrum of BPS branes.  This
concept reduces to variation of $\mu$-stability for B branes in the
large volume limit and based on this analogy, the work of Joyce on
such transitions for A branes \cite{Joyce}, and many other
considerations, a deformation of stability called $\Pi$-stability was
proposed in \cite{DFR1} which governs the spectrum of BPS branes on
arbitrary CY's in string theory.  This proposal has undergone a number
of tests; we also briefly describe work in progress which attempts to
complete the proposal by associating natural $t$-structures to regions
in K\"ahler moduli space.

Another front on which physical considerations have made progress is in
the explicit construction of these BPS branes and their moduli spaces.
A particularly simple class of constructions are the orbifold
constructions.  These can be exactly solved for orbifolds of flat
space, and provide a physical context which contains the structures
found both in orbifold resolution by quiver varieties \cite{Kron} and
in the study of the generalized McKay correspondence 
\cite{BKR,IN,Reid},
including generalizations of Beilinson's construction of
holomorphic bundles on $\BP^n$. \cite{Beil}

More recently, it has been found that results from the study of
explicit Gepner model boundary states \cite{RS} can be rederived starting
with a ``Landau-Ginzburg orbifold theory.''  \cite{DD}
Geometrically, this uses
an embedding of the CY in a weighted projective
space, and describes bundles on CY by
using the generalized McKay approach to describe bundles
on the ambient space, and then restricting.  One new point which emerges
is that the physical construction naturally contains moduli which
appear only on restriction.  Another striking result from this approach
is an explicit prediction for the connection formula which
relates the natural bases of
periods at the large volume and Gepner points.

\section{Background on CFT on Calabi-Yau}

The most general definition of a Dirichlet brane in weakly coupled
type \II\ string theory is that it is a conformal boundary condition
in the world-sheet CFT.  This definition is
discussed at this level of generality in \cite{Pol}.  It is quite
analogous to the discussion of branes in a geometrical framework as
calibrated submanifolds carrying appropriate vector bundles in the
sense that they are auxiliary objects; one requires some description
of the ambient space (the CY for us) to get started on their
definition.  

In the most general case, where we take as the ``ambient space'' an
abstract CFT defined in terms of a Hilbert space, Hamiltonian and
operator product coefficients, our starting point is a priori
non-geometric: it does not come with any local coordinates or other
conventional geometric data.  On general grounds, the CFT's which are
easy to describe in this abstract way are ``highly stringy'' -- the
conventional equations of physics such as the Einstein and Yang-Mills
equations will not be valid (except in special cases with high
symmetry) but will be drastically modified in string theory (we will
say more about this below).  In some cases, such as the Gepner models,
we can argue that they are still connected to solutions of Einstein's
equations by varying parameters (this is the ``large volume limit''),
and the interesting question is then
what aspects of this conventional geometrical interpretation survive
the stringy modifications and how do other aspects change.  However,
there are other cases such as asymmetric orbifolds \cite{NSV} where we
have no such large volume limit or other geometrical interpretation at
present.  In these cases, the abstract definition of D-brane makes
sense, but the identification of D-branes with cycles carrying bundles
on some space cannot even get started.

Although general, these remarks are intended to make the point that
many of the statements about branes which are taken as definition or
otherwise manifest in the mathematical treatments, actually require
proof or at least justification in the physical discussion starting
from string theory.  While we don't want to belabor this point, it
should be kept in mind.

We now specialize to the particular case of string compactification on
a Calabi-Yau threefold, for which the CFT is a $(2,2)$
superconformal field theory with $\hat c=3$.  These have been
discussed in
\cite{CK,Dijk,Morr,Voisin} and elsewhere.  The notation $(2,2)$
means there are two independent $N=2$ superconformal algebras describing
left and right-moving excitations on the string world-sheet.  Each $N=2$
algebra has two supercharges $G^+$ and $G^-$ which anticommute to produce
the Hamiltonian and a $U(1)$ charge $J$, which provides a grading on the
Hilbert space.

A distinguishing feature of the $\CN=2$ algebra is that it admits
a topological twisting, a redefinition of the generators after which one
of the supercharges can be used as a ``BRST charge,'' a differential $Q$
which squares to zero.  The cohomology of the sum $Q_L+Q_R$ of these
operators in the two $N=2$ algebras is finite dimensional and is the
Hilbert space of the topologically twisted theory.  This Hilbert space
also parameterizes a space of ``operators'' in the theory which correspond
to linearized deformations of the theory; these always turn out to 
integrable in the $(2,2)$ case so these theories come in families
with locally smooth moduli spaces.

Physicists are particularly interested in these models not because they
admit topological twisting but because they lead to string theories with
space-time supersymmetry.  Seeing this requires discussing a bit more
structure, namely the ``bosonized $U(1)$.''  The
main point here is that there are natural ``spectral flow'' operators denoted 
$e^{i\hat c\phi/2}$
(on left and right) which provide automorphisms between sectors of
different $U(1)$ charge $q$ and $q+\hat c/2$.  Sectors with integral
(resp. half-integral) $U(1)$ charge correspond to space-time bosons 
(fermions) and this one-to-one relation guarantees space-time supersymmetry.
The topologically twisted world-sheet theories are important in this context
mainly because their observables correspond to certain distinguished
amplitudes in this space-time supersymmetric theory.

The case of threefolds
with $\hat c=3$ is further distinguished in that the correlator of
three operators (the ``Yukawa coupling'')
is always well-defined; one can further show that such correlators
are the third derivatives of a function on moduli space $F$
known in the physics literature
as the ``$\CN=2$ prepotential'' and in the math literature as the
``Gromov-Witten potential.''  We will discuss the physical interpretation
of this potential below.

One made a choice in this construction of whether to identify $G^+$ or
$G^-$ as the BRST operator $Q$; it turns out that only the relative
choice between the two $N=2$ algebras matters in the final result.
This choice distinguishes ``A'' and ``B'' topologically twisted
theories, which although they arise from the same CFT, need have
little else in common.  On the level of abstract CFT, one can also
exchange $G^+$ and $G^-$ in one of the $\CN=2$ algebras to define a
``new'' CFT, all of whose physical predictions would be the same, but
in which the A and B twisted theories would be exchanged.  Thus there
is no a priori difference between A and B, and this is the sense in
which mirror symmetry is manifest in CFT.

This formalism arises from several more concrete constructions.
Contact with geometry and the large volume limit is most direct in
what is called the ``nonlinear sigma model'' definition.  This starts
from a Calabi-Yau threefold with Ricci-flat metric $g^{(0)}$,
and maps of the string world-sheet into this ``target space.''
The Ricci-flat metric is determined by a choice of complex structure (let
the moduli space of these be $MC(\CM)$) and complexified K\"ahler class
(moduli space $MK_0(\CM)$); physical arguments strongly suggest that
this produces a $(2,2)$ CFT determined by these moduli.
The two gradings provided by left and right $U(1)$ charge essentially
correspond to the bigrading of Dolbeault cohomology, and the B and A
twisted theories as above correspond respectively to a theory
describing variation of Hodge structure (the Yukawa couplings are the
second derivatives of the periods of the holomorphic three-form) and to
the ``quantum cohomology theory'' whose prepotential generates the
Gromov-Witten invariants.

Mirror symmetry equates the A model on $\CM$, which sees only $MK(\CM)$,
with the B model on its mirror $\CW$, which only sees $MC(\CW)$.
In particular, the complexified K\"ahler cone $MK_0(\CM)$
is only a large volume approximation to the ``true'' stringy K\"ahler
moduli space $MK(\CM)$, which is best defined simply as $MC(\CW)$.  
Much work has been done on the physical
interpretation of this point.  \cite{Greene}
Since there is such a close connection between these
topologically twisted theories and the complex geometry, the claim that 
these originate from physical CFT's for which mirror symmetry is manifest
places this formulation of mirror symmetry essentially beyond
doubt; various explicit arguments have confirmed this.  \cite{Morr}

Beyond these results which can be convincingly justified in the
topologically twisted string theory, what one actually has in this
definition is a way to compute general observables as a series
expansion in powers of the curvature of the metric $g^{(0)}$
multiplied by a distinguished scale of length $l_s$, the ``string
length.''  One also has prescriptions for adding ``instanton''
corrections associated to non-trivial holomorphic maps from
world-sheets of various genus into the CY.  It should be realized that
this expansion by itself does not provide a very convincing argument
for existence of these models on finite size CY, because such series
expansions often do not converge.  However, given evidence from other
directions that models do indeed exist, we can interpret them
as valid descriptions at least of the asymptotics of these observables
for large CY's.  (Whenever we talk about ``large'' CY's in the following,
we simply mean the case in which the leading terms are good approximations.)

Perhaps the most basic
observation one can make from these results is that
a CY is actually not Ricci flat in string theory  \cite{Grisaru}.  
One traditionally
derives Ricci-flatness as the condition for ``zero beta function''
required for conformal invariance; however this is just the leading term
and one sees the next correction at order $l_s^6$; the true condition
takes the form
\begin{equation}\label{beta-function}
0 = {\partial\over\partial \mu} g^\beta_{ij} =
R_{ij} + l_s^6 (R^{(4)})_{ij} + \ldots .
\end{equation}
with the term $R^{(4)}$ given in \cite{Grisaru}, 
and a higher order
term shown to be non-zero in \cite{JJP}.  It has also been shown that the
corrections are such that a solution can be found
at each order in the series \cite{NS}.

A conceptual argument that Ricci-flatness is not to be expected is simply
that the structure of the K\"ahler moduli space predicted by mirror
symmetry is quite different from the geometric K\"ahler cone or a naive
complexification of this, and we would expect the ``true'' stringy metric
to depend naturally on this true moduli space; the space of Ricci-flat
metrics does not.

However, it is difficult to go further with this definition, simply because
neither the equation (\ref{beta-function}) nor the metric $g^\beta$ it
describes are canonically defined; the theory of renormalization which
led to it does not distinguish in any convincing way
between $g^\beta$ and any other metric
$g^{\beta'}$ which would be produced by an arbitrary local functional
redefinition $g\rightarrow g + \alpha_1 R + \alpha_2 R^2 + \ldots$.

The study of D-branes will provide a better definition of the metric, 
which we discuss in the next section.  Let us turn now to alternatives
to the non-linear sigma model.  The most famous of these is the Gepner
model, also described in the references above.
This starts off from completely non-geometric data, but is believed
to be equivalent to the ``analytic continuation'' of the CY non-linear sigma
model to a specific highly stringy point in the true K\"ahler moduli space. 
The fact that such models do exist and are described by the $\CN=2$
formalism is our best reason to believe that the series expansions
described earlier can indeed be summed in some way to
define the models under discussion.

The best understanding we have at present of how the Gepner model is
connected to the geometric picture comes through Witten's ``linear
sigma model'' definition 
(\cite{Witten}; see \cite{GJS,HIV} for work on D-branes in this framework). 
One can think of this as defining the CY
through a specific embedding in $\BC^6$.  This embedding is described
in two steps.  First, one constructs a weighted projective space 
as a symplectic quotient of $\BC^6$ by a $U(1)$
action, in coordinates $z^i \rightarrow e^{iw_i\theta} z^i$ for chosen
weights $w_i$.  One then defines the CY as a complex hypersurface in
this ambient space.  
(One can similarly get more general toric varieties, etc.)

The symplectic quotient is the ordinary quotient by $U(1)$ of the
stable points in the preimage $\pi^{-1}(0)$ 
of a moment map $\pi:\BC^6\rightarrow\BR$
defined by
$$
\{z^i\} \rightarrow -\mu + \sum w_i |z^i|^2 .
$$
To get weighted projective space with weights $w_i$, one takes
$w_i\ge 0$ for $1\le i\le 5$, $w_6$ such that $\sum_i w_i=0$, and $\mu>0$.
Actually, the quotient is a line bundle $\CO_{\WP}(w_6)$ over the
weighted projective space; one then imposes the complex equations 
$\grad W =0$ with $W= z^6 f(z^i)$, which reduce
to $z^6=0$ and $f(z^i)=0$ within this zero section.

The important point is now that one can argue that varying the true
K\"ahler parameter to the Gepner point corresponds to an analytic
continuation to the ``Landau-Ginzburg phase,'' the limit $\mu<<0$.
Here the symplectic quotient is rather different: it is the orbifold
$\BC^5/\BZ_{-w_5}$, and the equations $\grad W=0$ now constrain the
coordinates to $z^i=0$ for $1\le i\le 5$.  Although on the surface this
looks quite different from the original nonlinear sigma model,
since the orbifold is birational to the
original line bundle, from the point of view of complex geometry the
two phases are not so different, and further physical considerations
bear this out, leading to a moduli space $MK$ covering both phases
with only isolated singularities, and no immediately evident distinction
between the phases.

This discussion motivates the idea that a good starting point for
defining CY is orbifold resolution, and indeed an even simpler class
of CY's can be defined this way, as the resolution of orbifolds $\BC^3/\Gamma$
with $\Gamma$ a discrete subgroup of $SU(3)$.  This produces non-compact
CY's, which have also been much studied as local models of singularities
in CY's, with the general conclusion being that mirror symmetry and
almost all of the qualitative features of CY physics are visible in these
simpler examples.  The simplest case is $\BC^3/\BZ_3$, studied in
\cite{DGM,DG,DFR2}.

\section{Dirichlet branes and stability}

The coarsest classification of boundary conditions in CFT is according
to the portion of the superconformal algebra they preserve.  One way
to specify a boundary condition is in terms of a linear transformation
implementing the ``reflection'' from left-moving to right-moving
excitations; this transformation must contain an automorphism of the
entire symmetry algebra and thus we need to know these.

In the $(2,2)$ context, there are two basic automorphisms one can use
to relate left and right moving algebras. \cite{OOY}\  The trivial choice
$G^{\pm}_L = G^\pm_R$ and $J_L=J_R$ is referred to as a ``B-type''
boundary condition as it is compatible with the B topological 
twisting.\footnote{This twisting is defined in open string conventions.}
Such boundary conditions can be obtained in the non-linear sigma model
by fixing the boundary to sit on a holomorphic submanifold of the CY
carrying a holomorphic bundle.  The other basic choice
is $G^{\pm}_L = G^{\mp}_R$ and $J_L=-J_R$ or ``A-type,'' which is
obtained by fixing the boundary to sit on a Lagrangian submanifold,
carrying a flat bundle.

This is all the data which is required to specify a boundary condition
in the topologically twisted open string theory \cite{Wittentop}, 
and in this context mirror symmetry
can be interpreted as a one-to-one equivalence between the two classes
of boundary conditions and all topological correlators.  This is the
physical background behind Kontsevich's proposal; however there is much
more to say about its physical underpinnings and meaning which has only
been addressed recently.

From a physical point of view, the most basic question is what
physical observables could one compute if one knew something about
the derived category of coherent sheaves or its mirror.
This question is not entirely trivial as the categorical framework 
is quite different from the usual physics language, which
assumes that the objects involved come in moduli spaces with explicit
coordinatizations, etc.

There are various physical realizations of the BPS branes we are
talking about, depending on what string theory we use and how many
flat dimensions the brane extends in.  All are essentially identical
in terms of the underlying string theory; the primary distinction is
between those which can be treated purely classically (so, objects are
solutions of equations of motion derived from some action), 
and those in which quantum effects
are taken into account (so objects are zero energy wave functions or
typically harmonic forms on the classical moduli space).  
In any case, they will
have half of the supersymmetry available in the original type \II\
theory on CY, because in that theory left and right-movers led to
independent supersymmetries, which are now related by the boundary
conditions.

Two realizations seem to be particularly useful to keep in mind.  In
the first, the physical picture is of open strings free to move in the
$4=10-6$ dimensions not taken up by the Calabi-Yau; this leads
to a world-volume theory with $\CN=1$ supersymmetry in $4$ dimensions,
which we treat classically.
These theories are gauge theories and the basic specification of such
a theory is a choice of gauge group (a semi-simple Lie group) $G$, a
K\"ahler manifold $X$ of ``chiral matter'' configurations admitting a
$G$-action by isometries, a $G$-invariant holomorphic function
$W:X\rightarrow \BC$ called ``superpotential,'' and finally moment
maps for the $G$ action.

Although all of this data is physical, for the primary question of
finding the BPS branes and their moduli space (as a complex variety),
the data only gets used in the following way: a configuration of the
brane system is a solution of $\grad W=0$ in the symplectic quotient
$X//G$.  We refer to the set of these as the moduli space of brane
configurations.
Since quotient commutes with restriction, one can also think
of this as the symplectic $G$-quotient of the variety $\grad W=0$.

The second realization takes the branes to be particles in the extra
$4$ dimensions.  This description is related to the first as follows:
given a moduli space of simple brane configurations (i.e. with endomorphism
group $U(1)$); one identifies particles with zero energy states in
a supersymmetric quantum mechanics on this moduli space.
By standard arguments, these will be cohomology
classes.  (If the moduli space is singular, it turns out that the
quotient/restriction definition gives it a natural embedding in a
non-singular space, which can be used to make the definition.)
These are the ``BPS branes'' of \cite{HM,GV} etc.

From now on we will consider only the first realization of brane, but use
the second to define what we will call ``BPS central charge'' and explain
the concept of marginal stability.  It can be shown \cite{Dijk} that
in an $\CN=2$ supersymmetric theory (such as type \II\ on CY),
the mass of a particle (in the usual space-time sense) satisfies a bound
\begin{eqnarray}\label{bps-bound}
M &\ge |Z(Q)|; \cr
Z(Q) &= q_i z^i + m^i {\partial F\over\partial z_i}
= \int_\Sigma \Omega
\end{eqnarray}
with equality attained only for BPS particles.  The quantity $Z(Q)$ is
the ``BPS central charge'' and depends only on the ``charge'' or
topological type of the brane and on the CY moduli.
It does not depend
on the particular point or cohomology class in the moduli space
of brane configurations.

For A branes, $Z(Q)$ is the integral of the holomorphic three-form over the
the Lagrangian submanifold, 
and thus depends only on the homology class of the submanifold
and the moduli $MC$.  The BPS bound is the same as that coming from
the calibrated geometry of special Lagrangian submanifolds. \cite{Harvey}

For B branes, $Z(Q)$ is defined in terms of the prepotential $F$, a function
on $MK(\CM)$, usually computed by invoking mirror symmetry, and depends
on the K-theory class of the brane and the moduli $MK$.  From now on, when
we speak of ``CY moduli space'' we mean the part ($MC$ or $MK$) which appears
in these considerations, and ``brane charge'' 
means the homology or K-theory class.

The most basic physical role of this prepotential is its appearance in
(\ref{bps-bound}).  Implicit in this setup is a flat (Gauss-Manin)
connection on the CY
moduli space, related to the choice of a fixed basis on homology.
This allows transporting a BPS brane between points in moduli space
and following the variation of its central charge.

Naively, this transport allows us to identify the complete spectrum
(list) of BPS branes at each point in moduli space.  However, this
is not true, because of the existence of
``lines of marginal stability'' on which the spectrum changes.
There is a strong constraint on these changes:
a brane with charge $Q$ can undergo a decay (and thus
disappear from the spectrum) into constituents
with charges $Q_i$ (satisfying $Q=\sum Q_i$) only if the phases
\begin{equation}\label{phase-def}
\varphi(Q_i) \equiv {1\over\pi}\Im\log Z(Q_i)
\end{equation}
are equal.  (The normalization factor $1/\pi$ will be explained later.)
This follows from energy conservation,
which implies that $M\ge \sum M_i$, and (\ref{bps-bound}), and
no other assumptions.  We will return to this point; for now we stress
that the spectrum of BPS branes does change with moduli, in both
realizations.  In the first realization, it will turn out that this
is the direct generalization of ``wall crossing'' phenomena related
to the behavior of $\mu$-stability under variations of K\"ahler class.

Once we know that moduli spaces of brane configurations exist, 
of course it would be
quite interesting to know more about them and in particular to find their
K\"ahler metrics.  In particular, the problem of how to give a canonical
definition of the metric of a CY in string theory has a natural answer
in this context: it is just the moduli space of the ``D0-brane,'' the
B brane which embeds in a point.  \cite{Doug1}
This presupposes that such a brane
exists and has a moduli space which is the CY, which although manifest
at large volume cannot
be taken for granted on a string-scale CY.  Furthermore, there might not be
a unique candidate; one can choose one by using the flat connection of
K\"ahler moduli space to carry the large volume D0-brane along some
path to the point of interest, a prescription which in general is path
dependent.  

Anyways, assuming that such a D0-brane exists, we stress that for each
point in the true K\"ahler moduli space (the complex structure moduli
space of the mirror), there exists a canonical K\"ahler metric on the
CY, which away from the large volume limit is not expected to be Ricci
flat (we will mention some explicit results bearing this out below
\cite{DG}).  This metric can also be thought of as the metric on
the moduli space of the D3-brane wrapped on the $T^3$ fibration of SYZ
mirror symmetry; again from this point of view, except in special
limits, it has no reason to be Ricci flat.  (In the explicit computations
of \cite{SYZ}, this would come about after instanton corrections.)
It is not even clear that
it has the same K\"ahler class as the original metric; indeed once
one goes from $MK_0(M)$ to $MK(M)$ one loses any clear sense of what
this would mean.

This metric or even the question of what equations determine it seems
rather inaccessible at present and we return to the question of what
Kontsevich's ideas should help physicists to compute.  In some sense,
the answer is the space $X$, group $G$ and superpotential $W$; however
it has taken some time for physicists to recognize any of this data in
Kontsevich and Fukaya's rather abstract derived categories, and indeed
this relation has not been completely spelled out in the literature.
The most basic point is simply that brane configurations can naturally
be thought of as objects in a category whose morphisms are the
linearized variations of the gauge theory data described above
(tangent vectors to a point in $X$ modulo gauge directions) associated
to the system which is the direct sum of two brane subsystems.  The
superpotential enters because such variations typically correspond to
obstructed deformations; the non-trivial point seems to be that in
cases associated to Calabi-Yau threefolds this obstruction theory can
always be summarized in terms of gradients of a potential.  This point
has been recognized in numerous special cases -- for example, the
functional $W$ for holomorphic bundles is the holomorphic Chern-Simons
functional, and for holomorphic curves is the functional proposed in
\cite{DT,WittenSuper} -- but a mathematical argument as
simple and general as the physical argument based on $\CN=1$ supersymmetry
does not seem to have appeared in the literature.

Having seen what Kontsevich's framework can describe, we also can see
in this language what it does not describe: namely, the details of the
symplectic quotient which are determined by the specific moment maps.
In particular, not all holomorphic objects (solutions of $\grad W=0$)
correspond to points in this quotient; only stable objects do (in
the GIT sense, which depends on the specific moment maps).  Thus, the
question of characterizing stable BPS branes (in the physical sense)
could be answered by proposing the mathematical stability
condition which they satisfy.  One can show by physical arguments that
just as the derived category of coherent sheaves does not depend on
the K\"ahler structure of the CY, the moment maps
should only depend on the K\"ahler structure and
not on the complex structure. \cite{Doug3}

On large CY's, branes correspond to bundles satisfying the hermitian
Yang-Mills equations, so the stability condition is just the one given by the
Donaldson and Uhlenbeck-Yau theorems, which is slope or $\mu$-stability:
a bundle $E$ (or coherent sheaf) is $\mu$-stable if, for all subsheaves
$E'$, one has $\mu(E')<\mu(E)$ with the slope $\mu(E)=c_1(E)/{\rm rank}(E)$.
This condition depends on both K\"ahler and complex structure but in a
clearly separated way: $\mu$ depends on K\"ahler structure, and the
subobject relation on complex structure.  This will be the
prototype for our stringy stability condition.

A first step towards such a condition is to consider
a deformed hermitian Yang-Mills equation derived from string theory
by Marino et. al. \cite{MMMS} as the
condition for unbroken supersymmetry from the ``Born-Infeld action.''
Their equation is a deformation of the Yang-Mills action which involves higher
powers of the Yang-Mills curvature $F$, and takes the form
\begin{equation}\label{MMMS-equation}
0 = \mu_{d,\theta} \equiv \Im e^{i\theta} \Tr (\omega + i l_s^2 F)^d
\end{equation}
on a $d$-fold with K\"ahler form $\omega$.  Considerations of geometric
invariant theory in \cite{Leung,Thom}
show that such an equation will have a solution 
for ``$\mu_\theta(E)$-stable'' bundles; i.e. a stability condition with 
$\mu_\theta \equiv \int \mu_{d,\theta}$ playing the role of the slope.

Although very explicit, this equation only takes into account power-like
corrections in $l_s$, while it is known that further instanton corrections
are present.  Related to this, it is not clear which K\"ahler form one
should take for $\omega$ (the Ricci-flat or some stringy version.)
This point could also be addressed by D-brane considerations, but it is
not clear that this is the best approach to take as 
one would prefer to have a condition which depends not on 
$\omega$ but on the point in $MK(M)$, i.e. 
on the complex structure of the mirror.

One can identify the precise quantity which generalizes slope
by appealing to mirror symmetry to relate this question to work of Joyce
on invariants counting special Lagrangian manifolds. \cite{Joyce}
Among the many results of this work is an analysis of the local stability
of a special Lagrangian manifold under variations of the complex structure;
this is determined by an inequality involving (in physics language)
the phase of the BPS central charge (\ref{phase-def}).  This fits
very well with the known physical considerations on marginal stability,
and the sign of the inequality is new information from this point of view.

This brings us to the proposal of \cite{DFR1}, that the BPS branes
at a specific point in (true) K\"ahler moduli space are the ``$\Pi$-stable''
branes, i.e. each brane based on a holomorphic object $E$ such that for
all subobjects $E'$ one has $\varphi(E')<\varphi(E)$.  Compared to Joyce's
criterion, the main difference is to use morphisms and subobjects in the
definition.  There is also a conceptual difference: instead of working
with special Lagrangians, we take the stability condition back to the B
picture and apply it there.  This has the great practical advantage that
the definition of subobject is fairly clear in the B picture (though we
will have to say more about this below) while computing the morphisms in
the A picture is a harder and still not fully understood problem.

Besides the arguments we cited, in \cite{DFR1} the proposal is
compared with more explicit results in the large volume and in
orbifold limits (more below) and finds agreement.  It was studied in
the example of the $\BC^3/\BZ_3$ orbifold and its resolution in
\cite{DFR2}, about which more below, and seems to produce sensible
results there.  As we mentioned, on physical grounds it is impossible
for the spectrum to change except on lines of marginal stability as
defined earlier, and it is hard to come up with any competing proposal
which satisfies this constraint.  Having said this, the proposal has
two ill-defined points.  First, one must lift the phases from the
interval $[0,2)$ (which is what (\ref{phase-def}) gives us) to $\BR$.
The geometry underlying this (in the A picture) is explained in
\cite{Seidel}; in our work we have only postulated this lifting ad
hoc.  Second, it requires objects to live in an abelian category (as
do all GIT notions of stability), but the only universal category
associated to a Calabi-Yau seems to be the derived category used by
Kontsevich, which is not abelian.  The two points are connected and we
will come back to them below.

For completeness, we should say that (in our opinion) the assumption
that A branes are in fact special Lagrangian submanifolds is on less
firm a footing than most of the other elements of the picture.  The
basic problem is the one that we mentioned, that the appropriate
definition of the metric and K\"ahler form of the CY in string theory
has not in fact been settled, but is almost certainly not the Ricci
flat metric.  One possibility is that one still has
the special Lagrangian condition but with respect to a preferred non
Ricci flat K\"ahler form, possibly the D0-metric described above.  
Some evidence for this idea can be found in the work of Leung, Yau
and Zaslow \cite{LYZ} which shows that the special Lagrangian condition is
related to the MMMS equation (\ref{MMMS-equation}), in which the metric
could be derived from D-brane considerations, by a Fourier-Mukai
transform of the type which should describe T duality in string theory.
At present however this is just a guess, and more physical work is
needed to understand this point; for example it is clear that with
some work a series expansion analogous to (\ref{beta-function}) could
be found for the corrections to the special Lagrangian condition.

\section{Flow of gradings}

We now turn to a point which emerges quite clearly from physics and CFT
but does not seem to have appeared in the mathematical discussion, namely
a ``flow of gradings'' which is induced under variation of K\"ahler
moduli (in the B picture) or complex structure (in the A picture).

Consider the B branes for definiteness; in the large volume limit
each brane corresponds to some stable coherent sheaf 
and thus for each pair of branes $E$ and $F$ we have 
graded morphisms $\Ext^n(E,F)$.  One often denotes this space as
$\Hom(E,F[n])$, a notation which finds its justification in homological
algebra.

Suppose we now follow a general path in K\"ahler moduli space 
between two points $K$ and $L$.
The new claim is that a morphism of degree $n$ at the starting point $K$
will undergo ``flow of grading'' determined by the phases (\ref{phase-def}): 
its degree at the point $L$ will be
$$
n' = n+\varphi_L(F)-\varphi_K(F)-\varphi_L(E)+\varphi_K(E) 
$$
(This is the reason for the $1/\pi$ in (\ref{phase-def}), which first
appeared in \cite{PZ}.)
In particular, the gradings are $\BR$-valued.  This rule can be expressed
more simply by saying that the grading of objects varies with $\varphi$,
so that the starting $\Hom(E[\varphi_K(E)],F[n+\varphi_K(F)])$ turns into 
$\Hom(E[\varphi_L(E)],F[n+\varphi_L(F)])$.

In itself this is just a definition but it can be checked physically if
we have more than one definition of the underlying CFT.  As we discuss
in the next section, this check has been made between large volume and
orbifold points in various models.  The phases are known from mirror
symmetry results and flow by integers between these points, and
the gradings of morphisms agree with these predictions.

The argument for the flow from CFT is quite simple. \cite{Doug3} It
relies on the identification of a morphism with a ``winding string in
the bosonized $U(1)$,'' and has the following intuitive picture.  One
can think of each brane as having location $\varphi$ on a circle of
circumference $2$ (the bosonized $U(1)$) and a morphism as an open string
stretched between the pair of branes; grading corresponds to length,
so the flow is simply induced by motion of the branes on the circle.

Once we grant this point, we quickly see that no single abelian
category (such as the category of coherent sheaves) could possibly
describe the branes throughout K\"ahler moduli space.  This is because
a general path will cause gradings to flow below zero, and convert
objects from even to odd grading, both of which lead to violations of
the axioms of an abelian category.  However, the derived category and
its distinguished triangles still make sense with these flows.  This
is perhaps the fundamental reason why one is forced to the derived
category in these considerations.  Recent work of Seidel
and Thomas and of Horja \cite{ST,Horja} 
provides even more concrete motivation for this point;
they show that the monodromies associated to general paths
in true K\"ahler moduli space act naturally on the derived category,
not the category of coherent sheaves.

One does not have a notion of subobject in the derived category and
thus no GIT definition of stability can apply in this context.  The
most direct way out of this problem is to propose that each point in
K\"ahler moduli space comes with a preferred abelian category.  This
could be defined by a choice of $t$-structure \cite{BBD}, which might
be determined by the phases as well: one basically wants to keep all
objects whose phases lie in the interval $(-1,0]$ and which are not
involved in morphisms of negative degree.  This point is presently
under investigation.

\section{Orbifolds}

The most accessible physical definitions of Calabi-Yau manifolds are the
orbifolds $\BC^3/\Gamma$ and the Landau-Ginzburg orbifolds described
in the introduction.  As it turns out, the basic constructions
involved in defining D-branes in these theories are already fairly
standard in mathematics and thus we can be relatively brief.

For B branes, one clearly wants to take some sort of
$\Gamma$-equivariant bundles on $\BC^n$, and this can be implemented
in physical terms by starting with the gauge theory Lagrangian for
branes in $\BC^n$ and applying a quotient which acts simultaneously on
$\BC^n$ and in the gauge group.  \cite{DM}
This leads directly to the quiver
theories which appear in the constructions of Kronheimer and
subsequent work.  \cite{Kron}  In particular, B branes are objects in a
``McKay'' quiver category, whose primitive objects (nodes) correspond
to irreps of $\Gamma$, whose links correspond to terms in the tensor
product with the representation defining the $\Gamma$ action on
$\BC^n$, and with specified quadratic relations.

The physical application requires $\Gamma\subset SU(n)$ and such orbifold
singularities can often be resolved (always for $n\le 3$) to smooth
manifolds.  This resolution is visible as the moduli space of a 
D0-brane and one can even compute the metric (actually, a controlled
approximation in the small blow-up limit); for $\BC^3/\Gamma$ one obtains
explicit non-Ricci-flat metrics of the type discussed earlier. \cite{DG}

More generally,
the question arises of how B branes as quiver objects or
irreps of $\Gamma$ are related to the large volume definition of B branes
as coherent sheaves.
In recent mathematical work \cite{BKR,IN,Reid}, 
a generalization of the McKay correspondence
has been developed, which in particular provides a natural basis for
K theory with compact support (thus supported on the exceptional divisor)
labelled by irreps of $\Gamma$.  The natural conjecture is that this
is the physical correspondence we are looking for. \cite{DD}

This conjecture can be tested against results
from mirror symmetry which determine the Gauss-Manin connection on
K\"ahler moduli space and lead to a connection formula which
relates the orbifold ``charge'' basis to the large volume basis.  This
gives an explicit expression for the Chern character of the bundle
corresponding to each irrep, and this has been checked to agree with
the conjecture for $\BC^3/\BZ_3$.  \cite{DFR2}
A somewhat more intricate version
of this applies to the Landau-Ginzburg orbifolds; one uses the McKay
correspondence to get the natural basis for the K theory of the
weighted projective space, and then restricts this to the Calabi-Yau.
Again one can compare with results from mirror symmetry and find
agreement. \cite{DD}

These constructions are intimately related to and generalize
Beilinson's construction of quiver categories from sheaves on $\BP^n$ 
\cite{Beil},
and it is these categories of quiver representations which are compared
to the category of coherent sheaves in the test of ``flow of gradings''
mentioned above.  More generally, all this
provides a fairly concrete relation between the problems of
classifying coherent sheaves and of solving certain $\CN=1$
supersymmetric gauge theories.  Although neither problem is easy in
general, the second problem is not only more familiar to physicists
but is a more appropriate language for the subsequent considerations
than an explicit classification of sheaves would have been.  It would
be even more valuable if such constructions could be found for all
sheaves, not just those which arise on restriction from the ambient space.

The first step in this generalization already appears in the physical
theories -- it turns out that these include additional moduli
which can be shown to correspond to linearized deformations of bundles
which appear after restriction to the CY. \cite{DD2}
Presumably one could go beyond this linearized level by computing
exact superpotentials; this may be possible using physical methods.

It would be quite interesting to know more about the variation of
the spectrum of BPS branes with K\"ahler moduli.  Some simple
results of this type were found for $\BC^3/\BZ_3$ in \cite{DFR2};
for example it was shown that there are lines of marginal stability
arbitrarily close to the large volume limit, by constructing bundles
which were $\epsilon$ away from violating slope stability.

It would clearly be important to find a microscopic derivation of
the stability condition, either from geometry (see \cite{Tstable}
for work in this direction) or perhaps CFT or string field theory.

Our general conclusion has to be that Kontsevich's proposal, coming as
it did before the study of D-branes, has turned out to be remarkably
prescient.  Although it has taken physicists some time to catch up,
we are finally extending this picture to provide both more concrete
pictures of branes on Calabi-Yau manifolds and new concepts which should
be of interest to mathematicians.

\medskip

We would like to thank D.-E. Diaconescu, B. Fiol and C. R\"omelsberger
for enjoyable collaborations, and
F. Bogomolov, M. Kontsevich, M. Marino, D. Morrison,
A. Polishchuk, P. Seidel, and R. P. Thomas for valuable discussions.
We thank R. P. Thomas for comments on the manuscript.

This work was supported by the Clay Mathematics Institute and
by DOE grant DE-FG02-96ER40959.

\end{document}